\documentclass{article}
\usepackage{amssymb}
\usepackage{amsbsy}
\begin{document}
\centerline{\Large \bf } \vskip 6pt

\begin{center}{\Large \bf Orlicz mixed chord integrals}\end{center}

\vskip 6pt\begin{center} \centerline{Chang-Jian Zhao\footnote{Research is supported by
National Natural Science Foundation of China
(11371334, 10971205).}}
\centerline{\it Department of Mathematics, China Jiliang
University, Hangzhou 310018, P. R. China}\centerline{\it Email:
chjzhao@163.com~~ chjzhao@cjlu.edu.cn}
\end{center}

\vskip 10pt

\begin{center}
\begin{minipage}{12cm}
{\bf Abstract}~ We introduce a affine geometric quantity and call
it {\it Orlicz mixed chord integral}, which generalize the chord
integrals to Orlicz space. Minkoswki and Brunn-Minkowski
inequalities for the Orlicz mixed chord integrals are establish.
These new inequalities in special cases yield some isoperimetric
inequalities for the usual chord integrals. The related concepts
and inequalities of $L_{p}$-mixed chord integrals are also
derived.


\end{minipage}
\end{center}
\vskip 20pt

\noindent{\large \bf 1 ~Introducation}\vskip 10pt

The radial addition $K\widetilde{+}L$ of star sets $K$ and $L$ can
be defined by
$$\rho(K\widetilde{+}L,\cdot)=\rho(K,\cdot)+\rho(L,\cdot),$$
where the star sets means compact sets that is star-shaped at $o$
and contains $o$ and $\rho(K,\cdot)$ denotes the radial function
of star set $K$. The radial function is defined by
$$\rho(K,u)=\max\{c\geq 0: cu\in K\},\eqno(1.1)$$ for $u\in S^{n-1}$, where $S^{n-1}$ denotes the
surface of the unit ball centered at the origin. The initial study
of the radial addition can be found in [1, P. 235]. $K$ is called
a star body, if $\rho(K,\cdot)$ is positive and continuous, and
let ${\cal S}^{n}$ denote the set of star bodies. Radial addition
and volume are the basis and core of dual Brunn-Minkowski theory
(see, e.g., [2], [3], [4], [5], [6], [7], [8] and [9]). It is
important that the dual Brunn-Minkowski theory can count among its
successes the solution of the Busemann-Petty problem in [3], [11],
[12], [13], and [14]. Recently, it has turned to a study extending
from $L_p$-dual Brunn-Minkowski theory to Orlicz dual
Brunn-Minkowski theory. The dual Orlicz-Brunn-Minkowski theory and
its dual have attracted people's attention [15], [16], [17], [18],
[19], [20], [21], [22], [23], [24] and [25].

For $K\in {\cal S}^{n}$ and $u\in S^{n-1}$, the half chord of $K$
in the direction $u$, defined by.
$$d(K,u)=\frac{1}{2}(\rho(K,u)+\rho(K,-u)).$$
If there exist constants $\lambda>0$ such that $d(K,u)=\lambda
d(L,u),$ for all $u\in S^{n-1}$, then star bodies $K,L$ are said
to have similar chord (see Gardner [1] or Schneider [26]). Lu [27]
introduced the chord integral of star bodies: For $K\in {\cal
S}^{n}$ and $0\leq i<n$, the chord integral of $K$, denoted by
$B_{i}(K),$ defined by
$$B_{i}(K)=\frac{1}{n}\int_{S^{n-1}}d(K,u)^{n-i}dS(u).\eqno(1.2)$$
For $i=0$, $B_{i}(K)$ becomes the chord integral $B(K)$.

The main aim of the present article is to generalize the chord
integrals to Orlicz the space. We introduce a new affine geometric
quantity such it Orlicz mixed chord integrals. The fundamental
notions and conclusions of the chord integrals and related
isoperimetric inequalities for the chord integrals are extended to
an Orlicz setting. The new inequalities in special case yield the
$L_p$-dual Minkowski, and Brunn-Minkowski inequalities for the
$L_{p}$-mixed chord integrals.

In Section 3, we introduce a notion of Orlicz chord addition
$K\check{+}_{\phi}L$ of star bodies $K,L$, defined by
$$\phi\left(\frac{d(K,x)}{d(K\check{+}_{\phi}L,x)}
,\frac{d(L,x)}{d(K\check{+}_{\phi}L,x)}\right)=1.\eqno(1.3)$$ Here
$\phi\in \Phi_{2}$, the set of convex function
$\phi:[0,\infty)^{2}\rightarrow (0,\infty)$ that are decreasing in
each variable and satisfy $\phi(0,0)=\infty$ and
$\phi(\infty,1)=\phi(1,\infty)=1$. The particular instance of
interest corresponds to using (1.5) with
$\phi(x_{1},x_{2})=\phi_{1}(x_{1})+\varepsilon\phi_{2}(x_{2})$ for
$\varepsilon>0$ and some $\phi_{1},\phi_{2}\in\Phi$, where the
sets of convex functions $\phi_{1}, \phi_{2}:[0,\infty)\rightarrow
(0,\infty)$ that are decreasing and satisfy
$\phi_{1}(0)=\phi_{2}(0)=\infty$,
$\phi_{1}(\infty)=\phi_{2}(\infty)=0$ and
$\phi_{1}(1)=\phi_{2}(1)=1$ .

In accordance with the spirit of Aleksandrov [28], Fenchel and
Jensen [29] introduction of mixed quermassintegrals, and
introduction of Lutwak's [30] $L_p$-mixed quermassintegrals, we
are based on the study of first order variational of the chord
integrals. In Section 4, we prove that the first order Orlicz
variation of the chord integrals can be expressed as: For $K,L\in
{\cal S}^{n},$ $\phi_{1},\phi_{2}\in {\Phi}$, $0\leq i<n$ and
$\varepsilon>0$,
$$\frac{d}{d\varepsilon}\bigg|_{\varepsilon=0^{+}}B_{i}
(K\check{+}_{\phi}\varepsilon\cdot L)=\frac{n-i}{n\cdot(\phi_{1})
'_{r}(1)}\cdot B_{\phi_{2},i}(K,L)^{n},\eqno(1.4)$$ where $\phi'_{r}(1)$ denotes the
value of the right derivative of convex function $\phi$ at point
$1$. In this first
order variational equation (1.4), we find a new geometric
quantity. Based on this, we extract the required geometric
quantity, denotes by $B_{\phi,i}(K,L)$ and call as Orlicz mixed
chord integrals, defined by
$$B_{\phi_{2},i}(K,L)=\left(\frac{n\cdot(\phi_{1})
'_{r}(1)}{n-i}\cdot\frac{d}{d\varepsilon}\bigg|_{\varepsilon=0^{+}}B_{i}
(K\check{+}_{\phi}\varepsilon\cdot L)\right)^{1/n}.\eqno(1.5)$$ We
also show the new affine geometric quantity has an integral
representation.
$$B_{\phi,i}(K,L)=\frac{1}{n}\int_{S^{n-1}}\phi
\left(\frac{d(L,u)}{d(K,u)}\right)d(K,u)^{n-i}dS(u).\eqno(1.6)$$
In Section 5, as application, we establish an Orlicz Minkowski
inequality for the Orlicz mixed chord integrals: If $K,L\in {\cal
S}^{n}$, $0\leq i<n$ and $\phi\in \Phi$, then
$$B_{\phi,i}(K,L)\geq
{B}_{i}(K)\cdot\phi\left(\left(\frac{{B}_{i}(L)}{{B}_{i}(K)}\right)^{1/(n-i)}\right)
.\eqno(1.7)$$ If $\phi$ is strictly convex, equality holds if and
only if $K$ and $L$ have similar chord. In Section 6, we establish
an Orlicz Brunn-Minkowski inequality for the Orlicz chord addition
and the chord integrals. If $K,L\in{\cal S}^{n}$, $0\leq i<n$ and
$\phi\in \Phi_{2}$, then
$$1\geq\phi\left(\left(\frac{{B}_{i}(K)}
{{B}_{i}(K\check{+}_{\phi}L)}\right)^{1/(n-i)},\left(\frac{{B}_{i}(L)}{{B}_{i}(K\check{+}_
{\phi}L)}\right)^{1/(n-i)}\right).\eqno(1.8)$$ If $\phi$ is
strictly convex, equality holds if and only if $K$ and $L$ have
similar chord.

\vskip 10pt \noindent{\large \bf 2 ~Preliminaries}\vskip 10pt

The setting for this paper is $n$-dimensional Euclidean space
${\Bbb R}^{n}$. A body in ${\Bbb R}^{n}$ is a compact set equal to
the closure of its interior. For a compact set $K\subset {\Bbb
R}^{n}$, we write $V(K)$ for the ($n$-dimensional) Lebesgue
measure of $K$ and call this the volume of $K$. Associated with a
compact subset $K$ of ${\Bbb R}^n$, which is star-shaped with
respect to the origin and contains the origin, its radial function
is $\rho(K,\cdot): S^{n-1}\rightarrow [0,\infty),$ defined by
$$\rho(K,u)=\max\{\lambda\geq 0: \lambda u\in K\}.$$
Note that the class (star sets) is closed under unions,
intersection, and intersection with subspace. The radial function
is homogeneous of degree $-1$, that is (see e.g. [1]),
$$\rho(K,ru)=r^{-1}\rho(K,u),$$
for all $u\in {S}^{n-1}$ and $r>0$. Let $\tilde{\delta}$ denote
the radial Hausdorff metric, as follows, if $K, L\in {\cal
S}^{n}$, then
$$\tilde{\delta}(K,L)=|\rho(K,u)-\rho(L,u)|_{\infty}.$$ From the definition of the radial function, it follows
immediately that for $A\in GL(n)$ the radial function of the image
$AK=\{Ay: y\in K\}$ of $K$ is given by (see e.g. [26])
$$\rho(AK,x)=\rho(K, A^{-1}x),\eqno(2.1)$$
for all $x\in {\Bbb R}^{n}$.

For $K_{i}\in {\cal S}^{n}, i=1,\ldots,m$, define the real numbers
$R_{K_{i}}$ and $r_{K_{i}}$ by
$$R_{K_{i}}=\max_{u\in S^{n-1}}d(K_{i},u),~~ {\rm and}~~ r_{K_{i}}=\min_{u\in S^{n-1}}d(K_{i},u),\eqno(2.2)$$
obviously, $0<r_{K_{i}}<R_{K_{i}},$ for all $K_{i}\in {\cal
S}^{n}$, and writing $R=\max\{R_{K_{i}}\}$ and
$r=\min\{r_{K_{i}}\}$, where $i=1,\ldots,m.$

\vskip 8pt {\it 2.1~ Mixed chord integrals}\vskip 8pt

If $K_{1},\ldots,K_{n}\in {\cal S}^{n}$, the mixed chord integral
of $K_{1},\ldots,K_{n}$, denotes by $B(K_{1},\ldots,K_{n})$,
defined by (see [27])
$$B(K_{1},\ldots,K_{n})=\frac{1}{n}\int_{S^{n-1}}d(K_{1},u)\cdots d(K_{n},u)dS(u).$$
If $K_{1}=\cdots=K_{n-i}=K,$ $K_{n-i+1}=\cdots=K_{n}=L$, the mixed
chord integral $B(K_{1},\ldots,K_{n})$ is written as $B_{i}(K,L)$.
If $L=B$ ($B$ is the unit ball centered at the origin), the mixed
chord integral $B_{i}(K,L)=B_{i}(K,B)$ is written as $B_{i}(K)$
and call chord integral of $K$. Obviously, For $K\in {\cal S}^{n}$
and $0\leq i<n$, we have
$$B_{i}(K)=\frac{1}{n}\int_{S^{n-1}}d(K,u)^{n-i}dS(u).\eqno(2.3)$$
If $K_{1}=\cdots=K_{n-i-1}=K,$ $K_{n-i}=\cdots=K_{n-1}=B$ and
$K_{n}=L$, the mixed chord integral
$B(\underbrace{K,\ldots,K}_{n-i-1},\underbrace{B,\ldots,B}_{i},L)$
is written as $B_{i}(K,L)$ and call $i$-th mixed chord integral of
$K$ and $L$. For $K,L\in {\cal S}^{n}$ and $0\leq i<n$, it is easy
that
$$B_{i}(K,L)=\frac{1}{n}\int_{S^{n-1}}d(K,u)^{n-i-1}d(L,u)dS(u).\eqno(2.4)$$
This integral representation (2.4), together with the H\"{o}lder
inequality, immediately gives: The Minkowski inequality for the
$i$-th mixed chord integral. If $K,L\in {\cal S}^{n}$ and $0\leq
i<n$, then
$$B_{i}(K,L)^{n-i}\leq B_{i}(K)^{n-i-1}B_{i}(L),\eqno(2.5)$$
with equality if and only if $K$ and $L$ have similar chord.

\vskip 8pt {\it 2.2~ $L_p$-mixed chord integrals}\vskip 8pt

Putting $\phi(x_{1},x_{2})=x_{1}^{-p}+x_{2}^{-p}$ and $p\geq 1$ in
(1.5), the Orlicz chord addition $\check{+}_{\phi}$ becomes a new
addition $\check{+}_{p}$ in $L_p$-space, and call as $L_p$-chord
addition of star bodies $K$ and $L$.
$$d(K\check{+}_{p}L,u)^{-p}=d(K,u)^{-p}+d(L,u)^{-p},\eqno(2.6)$$
for $u\in S^{n-1}$. The following result follows immediately form
(2.6) with $p\geq 1$.
$$-\frac{np}{n-i}\lim_{\varepsilon\rightarrow 0^{+}}\frac{B_{i}(K\check{+}_{p}
\varepsilon\cdot L)-B_{i}(L)}{\varepsilon}=\frac{1}{n}\int_{S^{n-1}}d(K,u)^{n-i+p}d(L,u)^{-p}dS(u).$$

{\bf Definition 2.1}~ Let $K,L\in {\cal S}^{n}$, $0\leq i<n$ and
$p\geq 1$, the $L_p$-chord integral of star $K$ and $L$, denotes
by $B_{-p,i}(K,L)$, defined by
$$B_{-p,i}(K,L)=\frac{1}{n}\int_{S^{n-1}}d(K,u)^{n-i+p}d(L,u)^{-p}dS(u).\eqno(2.7)$$
Obviously, when $K=L$, the $L_p$-mixed chord integral
$B_{-p,i}(K,K)$ becomes the chord integral $B_{i}(K).$  This
integral representation (2.7), together with the H\"{o}lder
inequality, immediately gives:

{\bf Proposition 2.2}~ {\it If $K,L\in {\cal S}^{n}$, $0\leq i<n$
and $p\geq 1$, then
$$B_{-p,i}(K,L)^{n-i}\geq B_{i}(K)^{n-i+p}B_{i}(L)^{-p},\eqno(2.8)$$
with equality if and only if $K$ and $L$ have similar chord.}

{\bf Proposition 2.3}~ {\it If $K,L\in {\cal S}^{n}$, $0\leq i<n$
and $p\geq 1$, then
$$B_{i}(K\check{+}_{p}L)^{-p/(n-i)}\geq B_{i}(K)^{-p/(n-i)}+B_{i}(L)^{-p/(n-i)},\eqno(2.9)$$
with equality if and only if $K$ and $L$ are dilates.}

{\it Proof}~ From (2.6) and (2.7), it is easily seen that the
$L_p$-chord integrals is linear with respect to the $L_p$-chord
addition, and together with inequality (2.8) show that for $p\geq
1$
$$B_{-p,i}(Q, K\check{+}_{p}L)=B_{-p,i}(Q,K)+B_{-p,i}(Q,L)~~~~~~~~~~~~~~~~~~~~~~~~~~~~~~~~~~~~~~~~~$$
$$~~~~~~~~~~~~~~\geq B_{i}(Q)^{(n-i+p)/(n-i)}(B_{i}(K)^{-p/(n-i)}+B_{i}(L)^{-p/(n-i)}),$$
with equality if and only if $K$ and $L$ have similar chord.

Take $K\check{+}_{p}L$ for $Q$, recall that
$B_{p,i}(Q,Q)=B_{i}(Q)$, inequality (2.9) follows
easy.

\vskip 10pt \noindent{\large \bf 3 ~Orlicz chord addition}\vskip
10pt

Throughout the paper, the standard orthonormal basis for ${\Bbb
R}^{n}$ will be $\{e_{1},\ldots,e_{n}\}$. Let $\Phi_{n},$
$n\in{\Bbb N}$, denote the set of convex function
$\phi:[0,\infty)^{n} \rightarrow (0,\infty)$ that are strictly
decreasing in each variable and satisfy $\phi(0)=\infty$ and
$\phi(e_{j})=1$, $j=1,\ldots,n$. When $n=1$, we shall write $\Phi$
instead of $\Phi_{1}$. The left derivative and right derivative of
a real-valued function $f$ are denoted by $(f)'_{l}$ and
$(f)'_{r}$, respectively. We first define the Orlicz chord
addition.

{\bf Definition 3.1}~ Let $m\geq 2, \phi\in\Phi_{m}$, $K_{j}\in
{\cal S}^{n}$ and $j=1,\ldots,m$, define the Orlicz chord addition
of $K_{1},\ldots,K_{m}$, denotes by
$\check{+}_{\phi}(K_{1},\ldots,K_{m})$, defined by
$$d(\check{+}_{\phi}(K_{1},\ldots,K_{m}),u)=\sup\left\{\lambda>0:
\phi\left(\frac{d(K_{1},u)}{\lambda},\ldots,\frac{d(K_{m},u)}{\lambda}\right)\leq
1\right\},\eqno(3.1)$$ for $u\in S^{n-1}.$ Equivalently, the
Orlicz chord addition $\check{+}_{\phi}(K_{1},\ldots,K_{m})$ can
be defined implicitly by
$$\phi\left(\frac{d(K_{1},u)}{d(\check{+}_{\phi}(K_{1},
\ldots,K_{m}),u)},\ldots,\frac{d(K_{m},u)}{d(\check{+}_{\phi}(K_{1},
\ldots,K_{m}),u)}\right)=1,\eqno(3.2)$$ for all $u\in {S}^{n-1}$.

An important special case is obtained when
$$\phi(x_{1},\ldots,x_{m})=\sum_{j=1}^{m}\phi_{j}(x_{j}),$$
for some fixed $\phi_{j}\in \Phi$ such that
$\phi_{1}(1)=\cdots=\phi_{m}(1)=1$. We then write
$\check{+}_{\phi}(K_{1},\ldots,K_{m})=K_{1}\check{+}_{\phi}\cdots\check{+}_{\phi}K_{m}.$
This means that $K_{1}\check{+}_{\phi}\cdots\check{+}_{\phi}K_{m}$
is defined either by
$$d(K_{1}\check{+}_{\phi}\cdots\check{+}_{\phi}K_{m},u)=
\sup\left\{\lambda>0:\sum_{j=1}^{m}\phi_{j}\left(\frac{d(K_{j},u)}{\lambda}
\right)\leq 1\right\},\eqno(3.3)$$ for all $u\in {S}^{n-1}$, or by
the corresponding special case of (3.2).

{\bf Lemma 3.2}~ {\it The Orlicz chord addition $\check{+}_{\phi}:
({\cal S}^{n})^{m}\rightarrow {\cal S}^{n}$ is monotonic.}

{\it Proof}~  This follows immediately from
(3.1).

{\bf Lemma 3.3}~ {\it The Orlicz chord addition $\check{+}_{\phi}:
({\cal S}^{n})^{m}\rightarrow {\cal S}^{n}$ is $GL(n)$ covariant.}

{\it Proof}~  This follows immediately from (2.1) and (3.1).

This shows Orlicz chord addition $\check{+}_{\phi}$ is $GL(n)$
covariant.

{\bf Lemma 3.4}~ {\it Suppose $K,\ldots,K_{m}\in{\cal S}^{n}$. If
$\phi\in \Phi$, then
$$\phi\left(\frac{d(K_{1},u)}{t}\right)+\cdots+\phi\left(\frac{d(K_{m},u)}{t}\right)=1$$
if and only if}
$$d(\check{+}_{\phi}(K_{1},
\ldots,K_{m}),u)=t$$

{\it Proof}~  This follows immediately from definition
3.1.

{\bf Lemma 3.5}~ {\it Suppose $K_{m},\ldots,K_{m}\in{\cal S}^{n}$.
If $\phi\in \Phi$, then}
$$\frac{r}{\phi^{-1}(\frac{1}{m})}\leq d(\check{+}_{\phi}(K_{1},
\ldots,K_{m}),u)\leq\frac{R}{\phi^{-1}(\frac{1}{m})}.$$

{\it Proof}~  This follows immediately from Lemma 3.4.

{\bf Lemma 3.6}~ {\it The Orlicz chord addition $\check{+}_{\phi}:
({\cal S}^{n})^{m}\rightarrow {\cal S}^{n}$ is continuous.}

{\it Proof}~ This follows immediately from Lemma 3.5.

Next, we define the Orlicz chord linear combination on the case
$m=2$.

{\bf Definition 3.7}~ Orlicz chord linear combination
$\check{+}_{\phi}(K,L,\alpha,\beta)$ for $K,L\in {\cal S}^{n}$,
and $\alpha,\beta\geq 0$ (not both zero), defined by
$$\alpha\cdot\phi_{1}\left(\frac{d(K,u)}{d(\check{+}_{\phi}(K,L,\alpha,\beta),u)}\right)+
\beta\cdot\phi_{2}\left(\frac{d(L,u)}
{d(\check{+}_{\phi}(K,L,\alpha,\beta),u)}\right)=1,\eqno(3.4)$$
for all $u\in {S}^{n-1}$.

We shall write $K\check{+}_{\phi}\varepsilon\cdot L$ instead of
$\check{+}_{\phi}(K,L,1,\varepsilon)$, for $\varepsilon\geq 0$ and
assume throughout that this is defined by (3.1), if $\alpha=1,
\beta=\varepsilon$ and $\phi\in \Phi$. We shall write
$K\check{+}_{\phi}L$ instead of $\check{+}_{\phi}(K,L,1,1)$ and
call the Orlicz chord addition of $K$ and $L$.

\vskip 10pt \noindent{\large \bf 4 ~Orlicz mixed chord
integrals}\vskip 10pt

In order to define Orlicz mixed chord integrals, we need the
following Lemmas 4.1-4.4.

{\bf Lemma 4.1}~ {\it Let $\phi\in \Phi$ and $\varepsilon>0$. If
$K,L\in {\cal S}^{n}$, then} $K\check{+}_{\phi}\varepsilon\cdot
L\in {\cal S}^{n}.$

{\it Proof}~ This follows immediately from (3.4) and Lemma 3.5.

{\bf Lemma 4.2}~ {\it If $K,L\in {\cal S}^{n}$, $\varepsilon>0$
and $\phi\in \Phi$, then
$$K\check{+}_{\phi}\varepsilon\cdot L\rightarrow K\eqno(4.1)$$ as} $\varepsilon\rightarrow
0^{+}$.

{\it Proof}~ This follows immediately from (3.4)and noting that
$\phi_{2}$, $\phi_{1}^{-1}$ and $d$ are continuous
functions.

{\bf Lemma 4.3}~ {\it If $K,L\in {\cal S}^{n}$, $0\leq i<n$ and
$\phi_{1}, \phi_{2}\in \Phi$, then}
$$\frac{d}{d\varepsilon}\bigg|_{\varepsilon=0^{+}}d(K\check{+}_{\phi}\varepsilon\cdot
L,u)^{n-i}=\frac{n-i}{(\phi_{1})'_{r}(1)}\cdot\phi_{2}
\left(\frac{d(L,u)}{d(K,u)}\right)\cdot d(K,u)^{n-i}.\eqno(4.2)$$

{\it Proof}~ Form (3.4), (4.1), Lemma 4.2 and notice that
$\phi_{1}^{-1}$, $\phi_{2}$ are continuous functions, we obtain
for $0\leq i<n$
$$\frac{d}{d\varepsilon}\bigg|_{\varepsilon=0^{+}}d(K\check{+}_{\phi}\varepsilon\cdot
L,u)^{n-i}=\frac{n-i}{(\phi_{1})'_{r}(1)}\cdot\phi_{2}
\left(\frac{d(L,u)}{d(K,u)}\right)\cdot
d(K,u)^{n-i}.$$

{\bf Lemma 4.4}~ {\it If $\phi\in \Phi_{2}$, $0\leq i<n$ and
$K,L\in {\cal S}^{n}$, then}
$$\frac{(\phi_{1})'_{r}(1)}{n-i}\cdot\frac{d}{d\varepsilon}\bigg|_{\varepsilon=0^{+}}B_{i}
(K\check{+}_{\phi}\varepsilon\cdot L)
=\frac{1}{n}\int_{S^{n-1}}\phi_{2}
\left(\frac{d(L,u)}{d(K,u)}\right)\cdot d(K,u)^{n-i}
dS(u).\eqno(4.3)$$

{\it Proof}~ This follows immediately from (2.1) and Lemma
4.2.

Denoting by $B_{\phi,i}(K,L)$, for any $\phi\in\Phi$ and $0\leq
i<n$, the integral on the right-hand side of (4.4) with $\phi_{2}$
replaced by $\phi$, we see that either side of the equation (4.3)
is equal to $B_{\phi_{2},i}(K,L)$ and hence this new Orlicz mixed
chord integrals $B_{\phi,i}(K,L)$ has been born.

{\bf Definition 4.5}~ For $\phi\in \Phi$ and $0\leq i<n$, Orlicz
mixed chord integrals of star bodies $K$ and $L$,
$B_{\phi,i}(K,L)$, defined by
$$B_{\phi,i}(K,L)=:\frac{1}{n}\int_{S^{n-1}}\phi
\left(\frac{d(L,u)}{d(K,u)}\right)\cdot d(K,u)^{n-i}dS(u).
\eqno(4.4)$$

{\bf Lemma 4.6}~ {\it If $\phi_{1}, \phi_{2}\in \Phi$, $0\leq i<n$
and $K,L\in {\cal S}^{n}$, then}
$$B_{\phi_{2},i}(K,L)=\frac{(\phi_{1})'_{r}(1)}{n-i}\lim_{\varepsilon\rightarrow 0^+}
\frac{{B}_{i}(K\check{+}_{\phi}\varepsilon\cdot
L)-{B}_{i}(K)}{\varepsilon}.\eqno(4.5)$$

{\it Proof}~ This follows immediately from Lemma 4.4 and
(4.4).

{\bf Lemma 4.7} {\it If $K,L\in {\cal S}^{n}$, $\phi\in {\cal C}$
and any $A\in{\rm SL(n)}$, then for $\varepsilon>0$}
$$A(K\check{+}_{\phi}\varepsilon\cdot L)=(AK)\check{+}_{\phi}\varepsilon\cdot(AL).\eqno(4.6)$$

{\it Proof}~ This follows immediately from (2.1) and (3.1).

We easy find that $B_{\phi,i}(K,L)$ is invariant under
simultaneous unimodular centro-affine transformation.

{\bf Lemma 4.8}~ {\it If $\phi\in \Phi$, $0\leq i<n$ and $K,L\in
{\cal S}^{n}$, then for $A\in SL(n)$,}
$$B_{\phi,i}(AK,AL)=B_{\phi,i}(K,L).\eqno(4.7)$$

{\it Proof}~ This follows immediately from Lemmas 3.3 and 4.7.

\vskip 10pt \noindent{\large \bf 5 ~Orlicz chord Minkowski
inequality}\vskip 10pt

In this section, we need define a Borel measure in $S^{n-1}$,
denotes by $B_{n,i}(K,\upsilon),$ call as chord measure of star
body $K$.

{\bf Definition 5.1}~ Let $K\in {\cal S}^{n}$ and $0\leq i<n$, the
chord measure, denotes by $B_{n,i}(K,\upsilon),$ defined by
$$dB_{n,i}(K,\upsilon)=\frac{d(K,\upsilon)^{n-i}}{n{B}_{i}(K)}dS(\upsilon).
\eqno(5.1)$$

{\bf Lemma 5.2}~ (Jensen's inequality) {\it Let $\mu$ be a
probability measure on a space $X$ and $g: X\rightarrow I\subset
{\Bbb R}$ is a $\mu$-integrable function, where $I$ is a possibly
infinite interval. If $\phi: I\rightarrow {\Bbb R}$ is a convex
function, then
$$\int_{X}\phi(g(x))d\mu(x)\geq\phi\left(\int_{X}g(x)d\mu(x)\right).
\eqno(5.2)$$ If $\phi$ is strictly convex, equality holds if and
only if $g(x)$ is constant for $\mu$-almost all $x\in X$} (see
[31, p.165]).

{\bf Lemma 5.3}~ {\it Suppose that $\phi: [0,\infty)\rightarrow
(0,\infty)$ is decreasing and convex with $\phi(0)=\infty$. If
$K,L\in{\cal S}^{n}$ and $0\leq i<n$, then
$$\frac{1}{n{B}_{i}(K)}\int_{S^{n-1}}\phi
\left(\frac{d(L,u)}{d(K,u)}\right)d(K,u)^{n-i}dS(u)\geq
\phi\left(\left(\frac{{B}_{i}(L)}{{B}_{i}(K)}\right)^{1/(n-i)}\right)
.\eqno(5.3)$$ If $\phi$ is strictly convex, equality holds if and
only if $K$ and $L$ have similar chord.}

{\it Proof}~ This follows immediately from (2.4), (2.5), (5.1) and Jensen's inequality.

{\bf Theorem 5.4}~ (Orlicz chord Minkowski inequality) {\it If
$K,L\in {\cal S}^{n}$, $0\leq i<n$ and $\phi\in \Phi$, then
$$B_{\phi,i}(K,L)\geq
{B}_{i}(K)\phi\left(\left(\frac{{B}_{i}(L)}{{B}_{i}(K)}\right)^{1/(n-i)}\right)
.\eqno(5.4)$$ If $\phi$ is strictly convex, equality holds if and
only if $K$ and $L$ have similar chord.}

{\it Proof}~ This follows immediately from (4.5) and Lemma 5.3.

{\bf Corollary 5.5} {\it If $K,L\in {\cal S}^{n}$, $0\leq i<n$ and
$p\geq 1$, then
$$B_{-p,i}(K,L)^{n-i}\geq{B}_{i}(K)^{n-i+p}{B}_{i}(L)^{-p},\eqno(5.5)$$
with equality if and only if $K$ and $L$ are dilates.}

{\it Proof} This follows immediately from Theorem 5.4 with
$\phi_{1}(t)=\phi_{2}(t)=t^{-p}$ and $p\geq 1$.

Taking $i=0$ in (5.6), this yields $L_{p}$-Minkowski inequality is
following: If $K,L\in {\cal S}^{n}$ and $p\geq 1$, then
$$B_{-p}(K,L)^{n}\geq B(K)^{n+p}B(L)^{-p},$$
with equality if and only if $K$ and $L$ have similar chord.

{\bf Corollary 5.6} {\it Let $K,L\in {\cal M}\subset{\cal S}^{n}$,
$0\leq i<n$ and $\phi\in \Phi$, and if either
$$B_{\phi,i}(Q,K)=B_{\phi,i}(Q,L),~ {\it for~ all}~ Q\in{\cal M}\eqno(5.6)$$
or
$$\frac{B_{\phi,i}(K,Q)}{{B}_{i}(K)}=\frac{B_{\phi,i}(L,Q)}{{B}_{i}(L)},~ {\it for~ all}~ Q\in{\cal M},\eqno(5.7)$$
then} $K=L.$

When $\phi_{1}(t)=\phi_{2}(t)=t^{-p}$ and $p\geq 1$, Corollary 5.6
becomes the following result.

{\bf Corollary 5.7} {\it Let $K,L\in {\cal M}\subset{\cal S}^{n}$,
$0\leq i<n$ and $p\geq 1$, and if either
$$B_{-p,i}(K,Q)=B_{-p,i}(L,Q),~ {\it for~ all}~ Q\in{\cal M}$$
or
$$\frac{B_{-p,i}(K,Q)}{{B}_{i}(K)}=\frac{B_{-p,i}(L,Q)}{{B}_{i}(L)},~ {\it for~ all}~ Q\in{\cal M},$$
then} $K=L.$

\vskip 10pt \noindent{\large \bf 6 ~Orlicz chord Brunn-Minkowski
inequality}\vskip 10pt

{\bf Lemma 6.1}~ {\it If $K,L\in {\cal S}^{n}$, $0\leq i<n$, and
$\phi_{1}, \phi_{2}\in \Phi$, then}
$${B}_{i}(K\check{+}_{\phi}L)=B_{\phi_{1},i}(K\check{+}_{\phi}L, K)
+B_{\phi_{2},i}(K\check{+}_{\phi}L, L).\eqno(6.1)$$

{\it Proof}~ This follows immediately from (3.1), (3.4) and (4.5).

{\bf Theorem 6.2}~ (Orlicz chord Brunn-Minkowski inequality)~ {\it
If $K,L\in{\cal S}^{n}$, $0\leq i<n$ and $\phi\in \Phi_{2}$, then
$$1\geq\phi\left(\left(\frac{{B}_{i}(K)}
{{B}_{i}(K\check{+}_{\phi}L)}\right)^{1/(n-i)},\left(\frac{{B}_{i}(L)}{{B}_{i}(K{+}_
{\phi}L)}\right)^{1/(n-i)}\right).\eqno(6.2)$$ If $\phi$ is
strictly convex, equality holds if and only if $K$ and $L$ have
similar chord.}

{\it Proof}~ This follows immediately from (5.4) and Lemma 6.1.

{\bf Corollary 6.3} {\it If $K,L\in {\cal S}^{n}$, $0\leq i<n$ and
$p\geq 1$, then
$${B}_{i}(K\check{+}_{p}L)^{-p/(n-i)}\geq{B}_{i}(K)^{-p/(n-i)}+{B}_{i}(L)^{-p/(n-i)},\eqno(6.3)$$
with equality if and only if $K$ and $L$ have similar chord.}

{\it Proof} The result follows immediately from Theorem 6.2 with
$\phi(x_{1},x_{2})=x_{1}^{-p}+x_{2}^{-p}$ and $p\geq 1$.

Taking $i=0$ in (6.3), this yields the $L_{p}$-Brunn-Minkowski
inequality for the
chord integrals. If $K,L\in{\cal S}^{n}$ and
$p\geq 1$, then
$$B(K\check{+}_{p}L)^{-p/n}\geq B(K)^{-p/n}+B(L)^{-p/n},$$
with equality if and only if $K$ and $L$ have similar chord.


\begin{thebibliography}{zz}
\vskip 8pt {\small
\bibitem{1} R. J. Gardner, Geometric Tomography, Cambridge Univ. Press, New York, 1996.
\bibitem{2} G. Berck, Convexity of $L_{p}$-intersection bodies, {\it Adv. Math.,} {\bf 222} (2009), 920-936.
\bibitem{3} R. J. Gardner, A. Koldobsky, T. Schlumprecht, An analytic solution to the Busemann-Petty problem on sections of
convex bodies, {\it Ann. Math.,} {\bf 149} (1999), 691-703.
\bibitem{4} C. Haberl, $L_{p}$ intersection bodies, {\it Adv. Math.,} {\bf 217} (2008),
2599-2624.
\bibitem{5} C. Haberl, M. Ludwig, A characterization of $L_{p}$ intersection
bodies, {\it Int. Math. Res. Not.,} 2006, Art. ID 10548, 29 pp.
\bibitem{6} A. Koldobsky, Fourier analysis in convex geometry, Mathematical
Surveys and Monographs 116, American Mathematical Society,
Providence, RI, 2005.
\bibitem{7} M. Ludwig, Intersection bodies and valuations, {\it Amer. J.
Math.,} {\bf 128} (2006), 1409-1428.
\bibitem{8} E. Lutwak, Centroid bodies and dual mixed volumes, {\it Proc. London Math. Soc.}, {\bf 60} (1990), 365-391.
\bibitem{9} E. M. Werner, R\'{e}nyi divergence and $L_{p}$-affine surface area for
convex bodies, {\it Adv. Math.,} {\bf 230} (2012), 1040-1059.
\bibitem{10} E. Lutwak, Dual mixed volumes, {\it Pacific J. Math.}, {\bf 58} (1975), 531-538.
\bibitem{11} R. J. Gardner, A positive answer to the Busemann-Petty
problem in three dimensions, {\it Ann. Math.,} {\bf 140}(2)
(1994), 435-447.
\bibitem{12} R. J. Gardner, A. Koldobsky, T. Schlumprecht, An analytic
solution to the Busemann-Petty problem on sections of convex
bodies, {\it Ann. Math.,} {\bf 149} (1999), 691-703.
\bibitem{13} F. E. Schuster, Valuations and Busemann-Petty type problems, {\it Adv.
Math.,} {\bf 219} (2008), 344-368.
\bibitem{14} E. Lutwak, Intersection bodies and dual mixed volumes, {\it Adv. Math.,} {\bf 71} (1988), 232-261.
\bibitem{15} R. J. Gardner, D. Hug, W. Weil, The Orlicz-Brunn-Minkowski theory: a general framework, additions, and
inequalities, {\it J. Differential Geom.,} {\bf 97}(3) (2014),
427-476.
\bibitem{16} E. Lutwak, D. Yang, G. Zhang, Orlicz projection bodies, {\it Adv. Math.,} {\bf 223} (2010), 220-242.
\bibitem{17} E. Lutwak, D. Yang, G. Zhang, Orlicz centroid bodies, {\it J. Differential Geom.,} {\bf 84} (2010), 365-387.
\bibitem{18} D. Xi, H. Jin, G. Leng, The Orlicz Brunn-Minkwski inequality, {\it Adv.
Math.,} {\bf 260} (2014), 350-374.
\bibitem{19} B. He, Q. Huang, On the Orlicz Minkowski problem for polytopes, {\it Discrete Comput. Geom.,} {\bf 48} (2012), 281-297.
\bibitem{20} C. Haberl, E. Lutwak, D. Yang, G. Zhang, The even Orlicz Minkowski problem, {\it Adv. Math.,} {\bf 224} (2010),
2485-2510.
\bibitem{21} J. Li, D. Ma, Laplace transforms and valuations, {\it J. Func. Anal.,} {\bf 272} (2017), 738-758.
\bibitem{22} Y. Lin, Affine Orlicz P\'{o}lya-Szeg\"{o} principle for log-concave functions, {\it J. Func. Aanl.}, {\bf 273} (2017), 3295-3326.
\bibitem{23} C. J. Zhao, On the Orlicz-Brunn-Minkowski theory, {\it Balkan J.
Geom. Appl.}, {\bf 22} (2017), 98-121.
\bibitem{24} C. J. Zhao, Orlicz dual mixed volumes, {\it Results Math.}, {\bf 68} (2015),
93-104.
\bibitem{25} C. J. Zhao, Orlicz dual affine quermassintegrals, {\it Forum Math.}, {\bf 30} (4) (2018), 929-945.
\bibitem{27} R. Schneider, Convex Bodies: The Brunn-Minkowski Theory, Second Edition, Cambridge Univ. Press, 2014.
\bibitem{28} F, Lu, Mixed chord-integrals of star bodies, {\it J.
Korean Math. Soc.}, {\bf 47} (2010) (2), 277-288.
\bibitem{29} A. D. Aleksandrov, Zur Theorie der gemischten Volumina von konvexen
K\"{o}rpern, I: Verall-gemeinerung einiger Begriffe der Theorie
der konvexen K\"{o}rper, {\it Mat. Sbornik  N. S.} {\bf 2}, (1937)
947-972.
\bibitem{30} W. Fenchel, B. Jessen, Mengenfunktionen und konvexe K\"{o}rper,
{\it Danske Vid Selskab Mat-fys Medd}, {\bf 16} (1938), 1-31.
\bibitem{31} E. Lutwak, The Brunn-Minkowski-Firey theory I: Mixed volumes and the Minkowski problem, {\it J. Differential Geom.}, {\bf 38} (1993), 131-150.
\bibitem{32} J. Hoffmann-J$\phi$gensen, Probability With a View Toward
Statistics, Vol. I, Chapman and Hall, New York, 1994, 165-243. }











\end{thebibliography}
\end{document}